# Delay-Robustness of Linear Predictor Feedback Without Restriction on Delay Rate

**Iasson Karafyllis[*] and Miroslav Krstic[**]**

[*]Dept. of Environmental Eng., Technical University of Crete,
73100, Chania, Greece, email: ikarafyl@enveng.tuc.gr

[**]Dept. of Mechanical and Aerospace Eng., University of California,
San Diego, La Jolla, CA 92093-0411, U.S.A., email: krstic@ucsd.edu

**Abstract**

Robustness is established for the predictor feedback for linear time-invariant systems with respect to possibly time-varying perturbations of the input delay, with a constant nominal delay. Prior results have addressed qualitatively constant delay perturbations (robustness of stability in $L^2$ norm of actuator state) and delay perturbations with restricted rate of change (robustness of stability in $H^1$ norm of actuator state). The present work provides simple formulae that allow direct and accurate computation of the least upper bound of the magnitude of the delay perturbation for which exponential stability in supremum norm on the actuator state is preserved. While prior work has employed Lyapunov-Krasovskii functionals constructed via backstepping, the present work employs a particular form of small-gain analysis. Two cases are considered: the case of measurable (possibly discontinuous) perturbations and the case of constant perturbations.

## 1. Introduction

Linear predictor feedback has been used widely for the stabilization of linear time-invariant systems with constant input delays. Artstein in [1] was the first to provide a rigorous extension of the so-called Smith predictor (see [9] and the discussion therein). Many applications and extensions of the linear predictor feedback have appeared in the literature (see for instance [12,13,14,15]). More recently, research efforts have been focused on nonlinear extensions of predictor-based feedback for nonlinear systems with input delays (see [2,4,6,7,8,9,10]) and the implementation issues of linear predictor feedback (see [17,18,19] and references therein).

However, the study of robustness properties of the linear predictor feedback with respect to perturbations of the input delay are rather scarce. To the best of our knowledge, the first robustness study for perturbations of the input delay appeared in [8], where Lyapunov techniques were employed. An alternative delay-robustness result for constant delays was presented in Section 5.3 in [9]. The efforts were continued in [2], where Lyapunov functionals were proposed for the robustness study for time-varying delays and perturbations. The results in [2] showed that not only the magnitude but the rate of change of the delay perturbation may be important for the robustness analysis. The norm on the actuator state in which stability was studied was $L^2$ in [8] and $H^1$ in [9, Section 5.3] and [2].

In this work, we consider the system:

$$\begin{aligned}\dot{x}(t) &= Ax(t) + Bu(t - r - \varepsilon d(t)) \\ x(t) &\in \Re^n, u(t) \in \Re^m, d(t) \in [-1,1]\end{aligned}, \text{ for } t \geq 0, \text{ a.e.} \quad (1.1)$$

where $0 < \varepsilon \leq r$ are constants. The linear predictor feedback is based on the constant nominal value of the delay $r > 0$:

$$u(t) = k\exp(Ar)x(t) + k\int_t^{t+r} \exp(A(t + r - s))Bu(s - r)ds, \text{ for } t \geq 0 \quad (1.2)$$

where $k \in \Re^{m \times n}$ is a constant vector such that the matrix $(A + Bk)$ is Hurwitz. We show that, provided that $\varepsilon > 0$ is sufficiently small, there exist constants $Q, \sigma > 0$ such that for all $x_0 \in \Re^n$, $u_0 \in C^0([-r - \varepsilon, 0]; \Re^m)$ with $u_0(0) = k\exp(Ar)x_0 + k\int_{-r}^0 \exp(-As)Bu_0(s)ds$ the solution $(x(t), u(t)) \in \Re^n \times \Re^m$ of (1.1), (1.2) with initial condition $x(0) = x_0$, $u(t) = u_0(t)$ for $t \in [-r - \varepsilon, 0]$ satisfies the following exponential stability estimate in the supremum norm of the actuator state:

$$|x(t)| + \max_{t - \varepsilon - r \leq s \leq t}(|u(s)|) \leq Q\exp(-\sigma t)\left(|x_0| + \max_{-\varepsilon - r \leq s \leq 0}(|u_0(s)|)\right), \forall t \geq 0 \quad (1.3)$$

for arbitrary disturbance $d : \Re_+ \to [-1,1]$ that belongs to one of the following classes:

1) The perturbation $d : \Re_+ \to [-1,1]$ is an arbitrary measurable function, i.e., $d \in L^\infty(\Re_+; [-1,1])$ (Theorem 2.1).

2) The perturbation $d : \Re_+ \to [-1,1]$ is constant (Corollary 2.3).

Clearly, (1.3) shows robust global exponential stability for the closed-loop system (1.1), (1.2). The estimation of $\varepsilon > 0$ will be given by explicit inequalities, which are derived by small-gain arguments. The inequalities can be used easily by the control practitioner in order to guarantee the successful application of the linear predictor feedback control strategy.

*Notation.* Throughout the paper we adopt the following notation:

∗ For a vector $x \in \Re^n$ we denote by $|x|$ its usual Euclidean norm, by $x'$ its transpose. For a real matrix $A \in \Re^{n \times m}$, $A' \in \Re^{m \times n}$ denotes its transpose and $|A| := \sup\{|Ax|; x \in \Re^n, |x| = 1\}$ is its induced norm. $I \in \Re^{n \times n}$ denotes the identity matrix.

∗ $\Re_+$ denotes the set of non-negative real numbers.

∗ Let $I \subseteq \Re$ be an interval and $U \subseteq \Re^m$ be a set. By $L^\infty(I; U)$ we denote the space of measurable and essentially bounded functions $u(\cdot)$ defined on $I$ and taking values in $U \subseteq \Re^m$. By $\sup_{t \in I}|u(t)|$ we denote the essential supremum for a function $u \in L^\infty(I; U)$. For $u \in L^\infty(I; U)$ we denote by $\|u\|$ the essential supremum of $u$ on $I \subseteq \Re$. If $I \subseteq \Re$ is an unbounded interval, then by $L^\infty_{loc}(I; U)$ we denote the space of measurable and locally essentially bounded functions $u(\cdot)$ defined on $I$ and taking values in $U \subseteq \Re^m$. By $C^0(I; U)$ we denote the space of continuous functions $u(\cdot)$ defined on $I$ and taking values in $U \subseteq \Re^m$.

## 2. Main Results

Arbitrary measurable perturbations $d \in L^\infty(\Re_+;[-1,1])$ of the delay can be considered for system (1.1). Indeed, we notice that this fact follows from the consideration of system (1.1) with

$$\dot{u}(t) = k\exp(Ar)(Ax(t) + Bu(t-r-\varepsilon d(t)) - Bu(t-r)) + kA\int_{-r}^{0}\exp(-As)Bu(t+s)ds + kBu(t) \quad (2.1)$$

Differential equation (2.1) is obtained by formally differentiating (1.2) with respect to $t \geq 0$. System (1.1) with (2.1) is a linear autonomous system described by Retarded Functional Differential Equations with disturbance $d \in L^\infty(\Re_+;[-1,1])$ and state space $\Re^n \times C^0([-r-\varepsilon,0];\Re^m)$ and satisfies all hypotheses (S1), (S2), (S3), (S4) in [5] for existence and uniqueness of solutions, for robustness of the equilibrium point and for the "Boundedness-Implies-Continuation" property. If we define the subspace

$$S := \left\{ (x,u) \in \Re^n \times C^0([-r-\varepsilon,0];\Re^m) : u(0) = k\exp(Ar)x + k\int_{-r}^{0}\exp(-As)Bu(s)ds \right\} \quad (2.2)$$

then we are in a position to guarantee that $S$ is a positively invariant set for system (1.1) with (2.1). Moreover, every solution of (1.1) with (2.1) and initial condition $(x_0, u_0) \in S$ is a solution of (1.1), (1.2) and every solution of (1.1), (1.2) with initial condition $(x_0, u_0) \in S$ is a solution of (1.1) with (2.1). Finally, we notice that there exist constants $M, L > 0$ such that for every $\varepsilon > 0$, $x_0 \in \Re^n$, $u_0 \in C^0([-r-\varepsilon,0];\Re^m)$, $d \in L^\infty(\Re_+;[-1,1])$ with $u_0(0) = k\exp(Ar)x_0 + k\int_{-r}^{0}\exp(-As)Bu_0(s)ds$ the unique solution $x \in C^0(\Re_+;\Re^n)$, $u \in C^0([-r-\varepsilon,+\infty);\Re^m)$ of system (1.1), (1.2) with initial conditions $x(0) = x_0$, $u(t) = u_0(t)$ for $t \in [-r-\varepsilon,0]$ satisfies the exponential growth estimate:

$$|x(t)| + |u(t)| \leq M\exp(Lt)\left(|x_0| + \max_{-r-\varepsilon \leq s \leq 0}|u_0(s)|\right), \quad \forall t \geq 0 \quad (2.3)$$

The existence of constants $M, L > 0$ satisfying estimate (2.3) follows directly from the integral representation of the solution of (1.1) with (2.1) and the Gronwall-Belman Lemma.

Our main result is the following theorem, which provides an explicit inequality for the magnitude $\varepsilon > 0$ of the delay perturbation under which robust global exponential stability for the closed-loop system (1.1), (1.2) is guaranteed.

**Theorem 2.1:** *Consider system (1.1), (1.2), where $0 < \varepsilon \leq r$ are constants, $A \in \Re^{n \times n}$, $B \in \Re^{n \times m}$, $k \in \Re^{m \times n}$ and $(A + Bk)$ is Hurwitz. There exist constants $Q, \sigma > 0$ such that for all $d \in L^\infty(\Re_+;[-1,1])$, $x_0 \in \Re^n$, $u_0 \in C^0([-r-\varepsilon,0];\Re^m)$ with $u_0(0) = k\exp(Ar)x_0 + k\int_{-r}^{0}\exp(-As)Bu_0(s)ds$ the solution $(x(t), u(t)) \in \Re^n \times \Re^m$ of (1.1), (1.2) with initial condition $x(0) = x_0$, $u(t) = u_0(t)$ for $t \in [-r-\varepsilon,0]$ satisfies estimate (1.3), provided that the following inequality holds:*

$$\Theta|\exp(Ar)Bk|\left(e^{|A+Bk|\varepsilon} - e^{-\lambda\varepsilon}\right) < \lambda \tag{2.4}$$

where $\Theta, \lambda > 0$ *are constants satisfying* $|\exp((A+Bk)t)| \leq \Theta e^{-\lambda t}$ *for all* $t \geq 0$. *Moreover, if* $n = 1$ *then inequality (2.4) can be replaced by the inequality*

$$2|Bk|\exp(Ar)\left(1 - e^{-|A+Bk|\varepsilon}\right) < |A+Bk| \tag{2.5}$$

**Remark 2.2:** Since the left hand-side of inequality (2.4) becomes zero for $\varepsilon = 0$, by continuity, there exists $\varepsilon > 0$ (sufficiently small) such that inequality (2.4) holds. The least upper bound value for $\varepsilon > 0$ can be determined numerically.

For the case of constant perturbations of the delay, we obtain the following result.

**Corollary 2.3:** *Consider the system*

$$\begin{aligned}\dot{x}(t) &= Ax(t) + Bu(t-\tau) \\ x(t) &\in \mathfrak{R}^n, u(t) \in \mathfrak{R}^m\end{aligned} \tag{2.6}$$

*with (1.2), where* $\tau, r \geq 0$ *are constants,* $A \in \mathfrak{R}^{n \times n}$, $B \in \mathfrak{R}^{n \times m}$, $k \in \mathfrak{R}^{m \times n}$ *and* $(A+Bk)$ *is Hurwitz. The zero solution of the closed-loop system is Globally Exponentially Stable provided that all roots of either of the following two equations:*

$$\det\left(sI - (A+Bk) + Bk\exp(Ar)\left(e^{-rs} - e^{-\varpi}\right)\right) = 0 \tag{2.7}$$

$$\det\left(sI - (A+Bk) + \exp(Ar)Bk\left(e^{-rs} - e^{-\varpi}\right)\right) = 0 \tag{2.8}$$

*have negative real parts.*

The following example illustrates the use of inequality (2.5) and Corollary 2.3.

**Example 2.4:** Consider the scalar system

$$\dot{x}(t) = x(t) + u(t-1-\varepsilon d(t)) \text{ with } x(t) \in \mathfrak{R}, u(t) \in \mathfrak{R}, d(t) \in [-1,1] \tag{2.9}$$

where $\varepsilon > 0$. For this example $A = 1 = B = r$ and we may choose $k = -p$, where $p > 1$. Theorem 2.1 guarantees that the closed-loop system (2.9) with

$$u(t) = -pe\, x(t) - p\int_0^1 e^s u(t-s)\,ds \tag{2.10}$$

and $d \in L^\infty(\mathfrak{R}_+; [-1,1])$ is robustly globally exponentially stable provided that $\varepsilon > 0$ satisfies

$$\varepsilon < \frac{1}{p-1}\ln\left(\frac{2pe}{2pe-p+1}\right) \tag{2.11}$$

In other words, system (2.9) with (2.10) is robustly globally exponentially stable provided that $\tau(t) \in (\tau_{\min}, \tau_{\max})$, where $\tau(t) = 1 + \varepsilon d(t)$, $\tau_{\min} = 1 - \varepsilon$, $\tau_{\max} = 1 + \varepsilon$ and $\varepsilon = \frac{1}{p-1}\ln\left(\frac{2pe}{2pe-p+1}\right)$.

On the other hand, if constant delay perturbations are considered, then the roots of the equation $s + (p-1) + pe^{1-\tau s} - pe^{1-s} = 0$ must have negative real parts. For every value of $p > 1$ there exist delay values $0 < \tau_{\min} < 1 < \tau_{\max}$ such that if $\tau \in (\tau_{\min}, \tau_{\max})$ then all roots of the equation $s + (p-1) + pe^{1-\tau s} - pe^{1-s} = 0$ have negative real parts. In order to determine the range of values of $\tau$ for which the roots of the equation $s + (p-1) + pe^{1-\tau s} - pe^{1-s} = 0$ have negative real parts, we determine the curves in the parameter plane (the $(p, \tau)$ plane) composed of points for which there exists $\omega \in \Re$ such that $\omega j + (p-1) + pe^{1-\tau \omega j} - pe^{1-\omega j} = 0$, where $j$ is the imaginary unit. The procedure that we follow for every $p > 1$, is:

(i) first we find numerically all solutions $\omega \in (0, 2pe)$ of the equation $(p-1)\cos(\omega) - \omega \sin(\omega) = \frac{(p-1)^2 + \omega^2}{2pe}$ (which is obtained from the equations $\cos(\omega\tau) - \cos(\omega) = -\frac{p-1}{pe}$ and $\sin(\omega\tau) - \sin(\omega) = \frac{\omega}{pe}$),

(ii) for every $\omega \in (0, 2pe)$ found from the previous step, we determine the unique solution $\phi \in \Re$ of the equations $\cos(\phi) = \cos(\omega) - \frac{p-1}{pe}$ and $\sin(\phi) = \sin(\omega) + \frac{\omega}{pe}$,

(iii) we find the positive solutions of $\tau = \frac{\phi + 2k\pi}{\omega}$, where $k$ is an arbitrary integer, and

(iv) finally, we collect all positive values of $\tau = \frac{\phi + 2k\pi}{\omega}$ from the previous step and we find the highest value of $\tau$ that is less than 1 (this is $\tau_{\min}$) and the lowest value of $\tau$ that is higher than 1 (this is $\tau_{\max}$).

The results are shown in Figure 1 both for measurable perturbations (where $\tau_{\min} = 1 - \varepsilon$, $\tau_{\max} = 1 + \varepsilon$ and $\varepsilon = \frac{1}{p-1}\ln\left(\frac{2pe}{2pe-p+1}\right)$) and for constant perturbations.

The bounds for the magnitude of the delay perturbation obtained from (2.11) are about 50% of the bounds obtained for constant perturbations. However, this is expected since (2.11) applies for measurable delay perturbations. Moreover, notice that the curves of $\tau_{\min}$ and $\tau_{\max}$ obtained for constant perturbations are not perfectly symmetric around 1. ◁

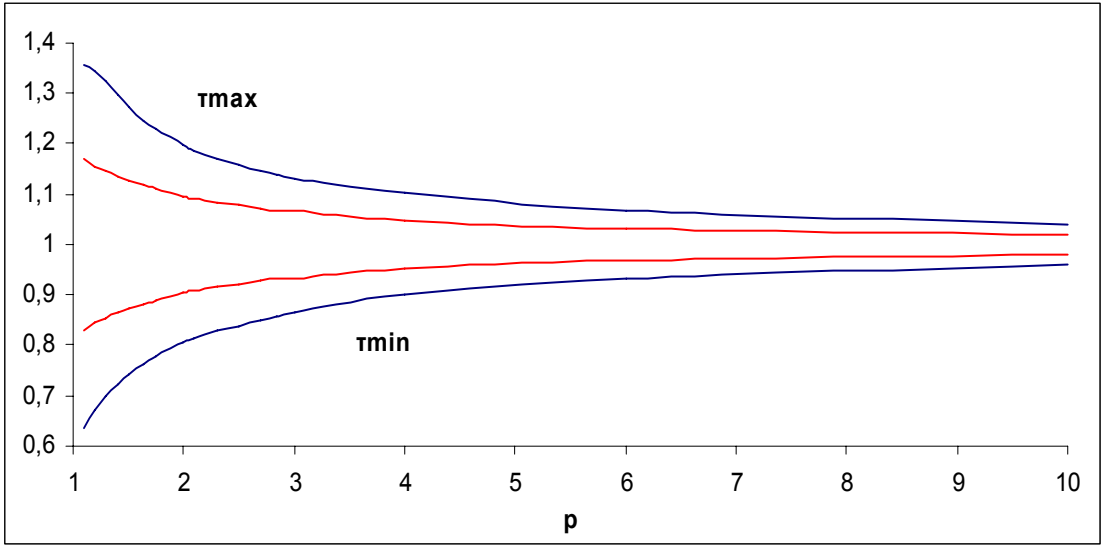

**Figure 1:** $\tau_{min}$ and $\tau_{max}$ for the closed-loop system (2.9) with (2.10). The red color is for measurable delay perturbations as calculated by (2.11) and the blue color for constant delay perturbations.

The proof of Theorem 2.1 relies on the following theorem.

**Theorem 2.5:** *Consider the system*

$$\begin{aligned}&\dot{x}(t) = Ax(t) + q(t)C\bigl(x(t-r-\varepsilon d(t)) - x(t-r)\bigr) \\ &x(t) \in \Re^n, d(t) \in [-1,1], q(t) \in [-1,1], u(t) \in \Re^n\end{aligned}, \text{ for } t \geq 0, \text{ a.e.} \qquad (2.12)$$

*where* $d \in L^\infty(\Re_+;[-1,1])$, $q \in L^\infty(\Re_+;[-1,1])$, $A,C \in \Re^{n \times n}$ *are constant matrices,* $r \geq \varepsilon \geq 0$ *are constants and* $A \in \Re^{n \times n}$ *is Hurwitz. Suppose that*

$$\Theta |C|\left(e^{|A|\varepsilon} - e^{-\lambda\varepsilon}\right) < \lambda \qquad (2.13)$$

*where* $\Theta, \lambda > 0$ *are constants satisfying* $|\exp(At)| \leq \Theta e^{-\lambda t}$ *for all* $t \geq 0$. *Then there exist constants* $Q, \sigma > 0$ *such that for all* $d \in L^\infty(\Re_+;[-1,1])$, $q \in L^\infty(\Re_+;[-1,1])$, $x_0 \in C^0([-r-\varepsilon,0];\Re^n)$ *the solution* $x(t) \in \Re^n$ *of (2.12) with initial condition* $x(t) = x_0(t)$ *for* $t \in [-r-\varepsilon,0]$ *that corresponds to inputs* $d \in L^\infty(\Re_+;[-1,1])$, $q \in L^\infty(\Re_+;[-1,1])$, *satisfies the following estimate*

$$|x(t)| \leq Q \exp(-\sigma t)\|x_0\|, \ \forall t \geq 0 \qquad (2.14)$$

*Moreover, if* $n=1$ *then inequality (2.13) can be replaced by the inequality*

$$2|C|\left(1 - e^{-|A|\varepsilon}\right) < |A| \qquad (2.15)$$

The proof of Theorem 2.5 is based on a small-gain argument and is provided in the following section. The small-gain argument for the proof of Theorem 2.5 was inspired by the results contained in [16], but the methodology of the proof is essentially different from that followed in [16].

**Remark 2.6:** The proof of Theorem 2.1 relies on showing the exponential stability properties of the system $\dot{p}(t) = (A + Bk)p(t) + \exp(Ar)Bk\big(p(t - r - \varepsilon d(t)) - p(t - r)\big)$, where $p(t) = \exp(Ar)x(t) + \int_{t}^{t+r} \exp(A(t + r - s))Bu(s - r)ds$ is the "predictor state". The exponential stability properties of the above system is guaranteed by means of Theorem 2.5. On the other hand Example 2.4 showed that the allowable magnitude of measurable delay perturbation is less than the magnitude obtained for constant perturbations from Corollary 2.3. We do not know if the conservatism is due to the small-gain approach (which is used for the proof of Theorem 2.5) or if the conservatism is due to the possibility that the stability analysis for delay perturbations depends not only on the magnitude of the perturbation but also on the rate of change of the perturbation. The latter implies that the rate of change of the perturbation may be important in stability analysis. Indeed, the recent work [2] has provided the construction of a Lyapunov functional for delay perturbations with constrained rate. Moreover, it should be noted that for delay perturbations with sufficiently small rate of change there exists a function $\phi : \Re_+ \to [0, r + \varepsilon]$, which satisfies $\phi(t) = r + \varepsilon d(t + \phi(t))$ for all $t \geq 0$: these are exactly the class of delays considered in [11] for which the following linear time-varying predictor feedback can be applied for the stabilization of (1.1):

$$u(t) = k\exp(A\phi(t))x(t) + k\int_{t}^{t+\phi(t)} \exp(A(t + \phi(t) - s))Bu(s - r - \varepsilon d(s))ds, \text{ for } t \geq 0 \tag{2.16}$$

provided that the function $d : \Re_+ \to [-1,1]$ is known.

## 3. Proofs of Main Results

We start with the proof of Theorem 2.5.

**Proof of Theorem 2.5:** If (2.13) holds, then (by continuity) there exists $\sigma \in (0, \lambda)$ such that:

$$e^{\sigma(r+\varepsilon)} \frac{\Theta|C|}{\lambda - \sigma}\Big(1 - e^{-(\lambda - \sigma)\varepsilon} + \big(\exp(|A|\varepsilon) - 1\big)\Big) < 1 \tag{3.1}$$

Let $d \in L^\infty(\Re_+;[-1,1])$, $q \in L^\infty(\Re_+;[-1,1])$, $x_0 \in C^0([-r-\varepsilon,0];\Re^n)$ be arbitrary and consider the solution $x(t) \in \Re^n$ of (2.12) with initial condition $x(t) = x_0(t)$ for $t \in [-r-\varepsilon, 0]$ that corresponds to inputs $d \in L^\infty(\Re_+;[-1,1])$, $q \in L^\infty(\Re_+;[-1,1])$. We define:

$$v(t) = x(t - r) - x(t - r - \varepsilon d(t)) \tag{3.2}$$

$$\|x\|_{[t_1,t_2]} := \max_{t_1 \leq s \leq t_2}\big(e^{\sigma s}|x(s)|\big), \quad \|v\|_{[t_1,t_2]} := \sup_{t_1 \leq s \leq t_2}\big(e^{\sigma s}|v(s)|\big) \tag{3.3}$$

for every $t_1 \leq t_2$ and we distinguish the following cases:

Case 1: $d(t) \leq 0$. In this case the following formula holds for the solution of system (2.12) for almost all $t \geq r$:

$$-v(t) = \left(\exp(A\varepsilon|d(t)|) - I\right)x(t-r) - \int_{t-r}^{t-r-\varepsilon d(t)} \exp(A(t-r-\varepsilon d(t)-s))q(s)Cv(s)ds \quad (3.4)$$

Using the fact that $|\exp(At)| \leq \Theta e^{-\lambda t}$ for all $t \geq 0$ and the fact that $|\exp(At) - I| \leq \exp(\|A\|\|t\|) - 1$, for all $t \in \Re$, we obtain from (3.4) for almost all $t \geq r$:

$$|v(t)|e^{\sigma t} \leq e^{\sigma r}\left(\exp(\|A\|\varepsilon) - 1\right)|x(t-r)|e^{\sigma(t-r)} + \Theta e^{\sigma r}\frac{1-e^{-(\lambda-\sigma)\varepsilon}}{\lambda-\sigma}|C| \sup_{t-r \leq s \leq t-r+\varepsilon}\left(e^{\sigma s}|v(s)|\right) \quad (3.5)$$

A direct consequence of definition (3.3) and inequality (3.5) is the following inequality which holds for all $t \geq r$:

$$\|v\|_{[r,t]} \leq e^{\sigma r}\left(\exp(\|A\|\varepsilon) - 1\right)\|x\|_{[0,t-r]} + \Theta e^{\sigma r}\frac{1-e^{-(\lambda-\sigma)\varepsilon}}{\lambda-\sigma}|C|\|v\|_{[0,t-r+\varepsilon]} \quad (3.6)$$

Case 2: $d(t) \geq 0$. In this case the following formula holds for the solution of system (3.13) for almost all $t \geq r + \varepsilon$:

$$v(t) = \left(\exp(A\varepsilon|d(t)|) - I\right)x(t-r-\varepsilon d(t)) - \int_{t-r-\varepsilon d(t)}^{t-r} \exp(A(t-r-s))q(s)Cv(s)ds \quad (3.7)$$

Similarly as in the previous case, using (3.7), we show that the following inequality holds for all $t \geq r + \varepsilon$:

$$\|v\|_{[r+\varepsilon,t]} \leq e^{\sigma(r+\varepsilon)}\left(\exp(\|A\|\varepsilon) - 1\right)\|x\|_{[0,t-r]} + \Theta e^{\sigma(r+\varepsilon)}\frac{1-e^{-(\lambda-\sigma)\varepsilon}}{\lambda-\sigma}|C|\|v\|_{[0,t-r]} \quad (3.8)$$

Consequently, we conclude from (3.6) and (3.8) that the following inequality holds for all $t \geq r + \varepsilon$:

$$\|v\|_{[r+\varepsilon,t]} \leq e^{\sigma(r+\varepsilon)}\left(\exp(\|A\|\varepsilon) - 1\right)\|x\|_{[0,t-r]} + \Theta e^{\sigma(r+\varepsilon)}\frac{1-e^{-(\lambda-\sigma)\varepsilon}}{\lambda-\sigma}|C|\|v\|_{[0,t-r+\varepsilon]} \quad (3.9)$$

Using the fact that $|\exp(At)| \le \Theta e^{-\lambda t}$ for all $t \ge 0$ and the variations of constants formula $x(t) = \exp(At)x(0) - \int_0^t \exp(A(t-s))q(s)Cv(s)ds$ for all $t \ge 0$, we obtain the estimate:

$$|x(t)|e^{\sigma t} \le \Theta e^{-(\lambda-\sigma)t}|x(0)| + \Theta \frac{1-e^{-(\lambda-\sigma)t}}{\lambda - \sigma}|C| \sup_{0 \le s \le t}\left(e^{\sigma s}|v(s)|\right), \text{ for all } t \ge 0 \qquad (3.10)$$

Definition (3.3) and inequality (3.10) in conjunction with the fact that $\sigma \in (0, \lambda)$ imply the following inequality:

$$\|x\|_{[0,t]} \le \Theta|x(0)| + \frac{\Theta|C|}{\lambda - \sigma}\|v\|_{[0,t]}, \text{ for all } t \ge 0 \qquad (3.11)$$

Combining (3.9) and (3.11), we obtain for all $t \ge r + \varepsilon$:

$$\|v\|_{[r+\varepsilon,t]} \le e^{\sigma(r+\varepsilon)}\left(\exp(|A|\varepsilon) - 1\right)\Theta|x(0)| + e^{\sigma(r+\varepsilon)}\frac{\Theta|C|}{\lambda - \sigma}\left(1 - e^{-(\lambda-\sigma)\varepsilon} + (\exp(|A|\varepsilon) - 1)\right)\|v\|_{[0,t]} \qquad (3.12)$$

Inequality (3.1) in conjunction with (3.12), implies the following inequality for all $t \ge 0$:

$$\|v\|_{[0,t]} \le e^{\sigma(r+\varepsilon)}\frac{\exp(|A|\varepsilon) - 1}{1 - \delta}\Theta|x(0)| + \|v\|_{[0,r+\varepsilon]} \qquad (3.13)$$

where $\delta := e^{\sigma(r+\varepsilon)}\frac{\Theta|C|}{\lambda - \sigma}\left(1 - e^{-(\lambda-\sigma)\varepsilon} + (\exp(|A|\varepsilon) - 1)\right) < 1$. Inequality (3.13) in conjunction with (3.11) and the fact that there exist constants $L, M > 0$ such that all solutions of (2.12) satisfy the estimate $|x(t)| \le M \exp(Lt) \max_{-r-\varepsilon \le s \le 0}|x(s)|$ and in conjunction with the fact that $\|v\|_{[0,r+\varepsilon]} \le 2e^{\sigma(r+\varepsilon)}\left(\|x\|_{[0,r+\varepsilon]} + \max_{-r-\varepsilon \le s \le 0}|x(s)|\right)$ (a direct consequence of definition (3.2)) imply that there exists a constant $Q > 0$ such that estimate (2.14) holds.

If $n = 1$ then $\Theta = 1$ and $\lambda = |A|$. If (2.15) holds then (by continuity) there exists $\sigma \in (0, |A|)$ such that $\delta := \frac{|C|e^{\sigma(r+\varepsilon)}}{|A| - \sigma}\left(2 - e^{-(|A|-\sigma)\varepsilon} - e^{-|A|\varepsilon}\right) < 1$. Moreover, inequalities (3.6) and (3.8) are replaced by the following inequalities:

$$\|v\|_{[r,t]} \le e^{\sigma r}\left(1 - \exp(-|A|\varepsilon)\right)\|x\|_{[0,t-r]} + e^{\sigma r}\frac{1 - e^{-(|A|-\sigma)\varepsilon}}{|A| - \sigma}|C|\|v\|_{[0,t-r+\varepsilon]}$$

$$\|v\|_{[r+\varepsilon,t]} \leq e^{\sigma(r+\varepsilon)}\left(1-\exp(-|A|\varepsilon)\right)\|x\|_{[0,t-r]} + e^{\sigma(r+\varepsilon)}\frac{1-e^{-(|A|-\sigma)\varepsilon}}{|A|-\sigma}|C|\|v\|_{[0,t-r]}$$

It follows that inequality (3.9) is replaced by

$$\|v\|_{[r+\varepsilon,t]} \leq e^{\sigma(r+\varepsilon)}\left(1-\exp(-|A|\varepsilon)\right)\|x\|_{[0,t-r]} + e^{\sigma(r+\varepsilon)}\frac{1-e^{-(|A|-\sigma)\varepsilon}}{|A|-\sigma}|C|\|v\|_{[0,t-r+\varepsilon]} \quad (3.14)$$

Combining (3.14) with (3.11) and $\Theta=1$, $\lambda=|A|$, we obtain the estimate:

$$\|v\|_{[r+\varepsilon,t]} \leq e^{\sigma(r+\varepsilon)}\left(1-e^{-|A|\varepsilon}\right)|x(0)| + \frac{|C|e^{\sigma(r+\varepsilon)}}{|A|-\sigma}\left(2-e^{-(|A|-\sigma)\varepsilon}-e^{-|A|\varepsilon}\right)\|v\|_{[0,t]}$$

Since $\delta := \frac{|C|e^{\sigma(r+\varepsilon)}}{|A|-\sigma}\left(2-e^{-(|A|-\sigma)\varepsilon}-e^{-|A|\varepsilon}\right)<1$, the above inequality implies the inequality $\|v\|_{[0,t]} \leq e^{\sigma(r+\varepsilon)}\frac{1-e^{-|A|\varepsilon}}{1-\delta}|x(0)| + \|v\|_{[0,r+\varepsilon]}$. The previous inequality in conjunction with (3.11) and the fact that there exist constants $L,M>0$ such that all solutions of (3.13) satisfy the estimate $|x(t)| \leq M\exp(Lt)\max_{-r-\varepsilon\leq s\leq 0}|x(s)|$ and in conjunction with the fact that $\|v\|_{[0,r+\varepsilon]} \leq 2e^{\sigma(r+\varepsilon)}\left(\|x\|_{[0,r+\varepsilon]} + \max_{-r-\varepsilon\leq s\leq 0}|x(s)|\right)$ (a direct consequence of definition (3.2)) imply that there exists a constant $Q>0$ such that estimate (2.14) holds. The proof is complete. ◁

We are now ready to provide the proof of Theorem 2.1.

**Proof of Theorem 2.1:** Let arbitrary $(x_0,u_0)\in S$ (where $S$ is defined by (2.2)), $d\in L^\infty(\Re_+;[-1,1])$ and consider the solution $(x(t),u(t))\in\Re^n\times\Re^m$ of (1.1), (1.2) with initial conditions $x(0)=x_0$, $u(t)=u_0(t)$ for $t\in[-r-\varepsilon,0]$ corresponding to $d\in L^\infty(\Re_+;[-1,1])$. Define for all $t\geq 0$:

$$p(t) = \exp(Ar)x(t) + \int_t^{t+r}\exp(A(t+r-s))Bu(s-r)ds \quad (3.15)$$

Notice that (1.2) and definition (3.15) implies that the following equality holds for all $t\geq 0$:

$$u(t) = kp(t), \text{ for all } t\geq 0 \quad (3.16)$$

By using (1.1) and definition (3.15), it follows that the following differential equation holds for all $t \geq 0$:

$$\dot{p}(t) = \exp(Ar)Ax(t) + \exp(Ar)Bu(t - r - \varepsilon d(t)) + A(p(t) - \exp(Ar)x(t)) \\ + Bu(t) - \exp(Ar)Bu(t - r) \tag{3.17}$$

Using the identity $A\exp(Ar) = \exp(Ar)A$ and (3.16) it follows that the following differential equation holds for all $t \geq r + \varepsilon$:

$$\dot{p}(t) = (A + Bk)p(t) + \exp(Ar)Bk(p(t - r - \varepsilon d(t)) - p(t - r)) \tag{3.18}$$

Inequalities (2.4) (or (2.5)) guarantee that Theorem 2.5 can be applied to system (3.18) and that there exist constants $Q, \sigma > 0$ such that the following inequality holds:

$$|p(t)| \leq Q\exp(-\sigma(t - r - \varepsilon)) \max_{0 \leq s \leq r + \varepsilon} |p(s)|, \quad \forall t \geq r + \varepsilon \tag{3.19}$$

Using (3.19) in conjunction with (3.15), (3.16), (2.3) and the following equality:

$$x(t) = \exp(-Ar)p(t) - \int_t^{t+r} \exp(A(t - s))Bkp(s - r)ds \tag{3.20}$$

which holds for all $t \geq r$ and is a direct consequence of (3.15) and (3.16), we obtain (1.3) (possibly with different constants $Q, \sigma > 0$). The proof is complete. ◁

Finally, we end this section with the proof of Corollary 2.3.

**Proof of Corollary 2.3:** The proof has two parts: in the first part we show that if all roots of equation (2.8) have negative real parts then the zero solution is Globally Exponentially Stable for system (2.6) with (1.2) and in the second part we show that if all roots of equation (2.7) have negative real parts then the zero solution is Globally Exponentially Stable for system (2.6) with (1.2).

<u>First part:</u> If all roots of equation (2.8) have negative real parts then Corollary 6.1 on page 215 in [3] guarantees that the zero solution is Globally Exponentially Stable for the system:

$$\dot{p}(t) = (A + Bk)p(t) + \exp(Ar)Bk(p(t - \tau) - p(t - r)) \tag{3.21}$$

Notice that the differential equation (3.21) holds for all $t \geq \max(r, \tau)$ for $p(t)$ as defined by (3.15) for all $t \geq 0$. Using (2.3), (3.16) and (3.20) we conclude that the zero solution is Globally Exponentially Stable for system (2.6) with (1.2).

Second part: Equation (1.2) implies for all $t \geq \tau$:

$$u(t-\tau) = k\exp(Ar)x(t-\tau) + k\int_{t-\tau}^{t-\tau+r}\exp(A(t-\tau+r-s))Bu(s-r)ds \qquad (3.22)$$

The variations of constants formula for (2.6) in conjunction with (3.22) implies that

$$u(t-\tau) = k\exp(Ar)x(t-\tau) + kx(t) - k\exp(Ar)x(t-r), \text{ for all } t \geq \max(r,\tau) \qquad (3.23)$$

Therefore, exponential stability for system for system (2.6) with (1.2) is guaranteed by the exponential stability of the system:

$$\dot{x}(t) = (A+Bk)x(t) + Bk\exp(Ar)(x(t-\tau) - x(t-r)) \qquad (3.24)$$

If all roots of equation (2.7) have negative real part then Corollary 6.1 on page 215 in [3] guarantees that zero solution is Globally Exponentially Stable for system (3.24). The proof is complete. ◁

## 4. Concluding Remarks

We have provided formulae that allow us to compute estimates of the least upper bound of the magnitude of the delay perturbation that does not destroy the exponential stability properties of the closed-loop system (1.1) with (1.2). Two cases have been considered: the case of measurable perturbations and the case of constant perturbations.

The formulae can be used easily by the control practitioner in order to estimate the delay error that can be tolerated. For the case of measurable time-varying perturbations the magnitude of the delay perturbation $\varepsilon > 0$ must satisfy inequality (2.4) (or (2.5) if $n=1$). All quantities involved in inequality (2.4) can be computed easily using software packages to compute the norms of matrices $\exp(Ar)Bk \in \Re^{n\times n}$, $(A+Bk) \in \Re^{n\times n}$ and to determine the constants $\Theta, \lambda > 0$ by finding a symmetric positive definite matrix $P \in \Re^{n\times n}$ and a constant $\mu > 0$ that satisfies $P(A+Bk) + (A+Bk)'P + 2\mu P \leq 0$ and $P \geq I$ (select $\lambda = \mu$ and $\Theta = \sqrt{|P|}$ ).

An example showed that the allowable magnitude of measurable delay perturbation is less than the magnitude obtained for constant perturbations from Corollary 2.3. We do not know if the conservatism is due to the small-gain approach (which is used for the proof of Theorem 2.5) or if the conservatism is due to the possibility that the stability analysis for delay perturbations depends not only on the magnitude of the perturbation but also on the rate of change of the perturbation. The latter implies that the rate of change of the perturbation may be important in stability analysis. It remains an open problem to construct more accurate expressions for the tolerance of the delay error which may involve the rate of change of the delay perturbation.

As in [8], where a Lyapunov analysis in $L^2$ is pursued, our stability analysis in $C^0$ separately considers positive and negative perturbations of the delay, whereas the Lyapunov analyses in $H^1$ in Section 5.3 in [9], and in [2] simultaneously tackle positive and negative perturbations on the delay.